\documentclass[a4paper]{amsart}

\usepackage[T1]{fontenc}
\usepackage[linktocpage]{hyperref} 
	\hypersetup{colorlinks=true,linkcolor=blue,citecolor=blue,filecolor=magenta,urlcolor=blue}      
\usepackage{geometry}       
\usepackage[english]{babel}  
\usepackage{amsmath,amsfonts,amssymb,amsthm}
\usepackage{tikz,tkz-graph,tikz-cd}

\newtheorem{thm}{Theorem}[section]

\theoremstyle{definition}

\theoremstyle{remark}
\newtheorem{rmk}[thm]{Remark}

\renewcommand{\epsilon}{\varepsilon}

\begin{document}

\title{Winding quotients for virtual period maps of rank 1}

\author[K. Saito]{Kyoji Saito}
			 \address{RIMS, Kyoto university}
			 \email{saito@kurims.kyoto-u.ac.jp}

\thanks{}				

\subjclass[2010]{
}		


\begin{abstract}
We illustrate a rank 1 model of  virtual period maps and their associated  winding quotients. The winding quotient is a new notion that appears in the recent study of virtual period maps due to the appearance of exponents which are imaginary numbers. 
It requires a reformulation of the classical inversion problem for the classical period maps. We answer to the new inversion problem by introducing a $q$-multiplicatively periodic function, whose pull-back to the winding covering space is the Weierstrass $\wp$-function up to a correction by Eisenstein series  $E_2$. This function appears as the propagator on elliptic curves in mathematical physics. 
\end{abstract}

\maketitle
\vspace{-0.3cm}

\tableofcontents

\vspace{-0.5cm}
\section{Introduction}
The classical theory of elliptic integral of the first kind for a Weierstrass family of elliptic curves $X\to S$ (e.g.\ \cite{Si} Chap.1)  defines a  period map, which is a multi-valent map from the complement of the discriminant divisor  $S\setminus D$ of the base space $S$ of the family to its range domain (so called, period domain).  
The range is isomorphic to the monodromy covering  of  $S\setminus D$, i.e.\  the quotient space of the universal covering  $(S\setminus D)^{\sim}$ divided by the kernel  of the monodromy representation of  $\pi_1(S\setminus D,*)$. Then, the Jacobi's inversion problem for the period map ([ibid] 1.13.) 
 is to ask for a description of the inverse map from the range to the defining domain $S\setminus D$, and its answer is given by Eisenstein series (see [ibid] Section13)

The theory of period integrals for primitive forms \cite{S1} has generalized those classical theory of  period integrals to  wider classes of geometric families $X\to S$. In particular, the elliptic integrals of Weierstrass type, Jacobi type and Hesse type are explained by the primitive forms associated to the root systems of types $\mathrm{A}_2, \mathrm{B}_2$ and $ \mathrm{G}_2$, respectively  (\cite{S2})). Recently, we find that  hyperbolic root systems of rank 2, which do not have geometric origin, admit also the construction of primitive forms and associated period maps (\cite{S3}). We call such period map ``virtual'' since it does not have corresponding geometric family $X\to S$. Then, in these virtual period maps, we observe some new phenomenon. Namely,  an infinite cyclic group, called the winding group, acts on the defining domain and on the monodromy covering space of the period map. It seems to  ask for a  reformulation of the classical inversion problem. 
 

The purpose of the present note is to illustrate such winding phenomenon by an elementary  model of virtual period maps of rank 1.  We observe that,  in this model,  the range is the quotient   of the monodromy covering space by the action of the winding group,  which we call the winding quotient range. Thus, the classical inversion problem is reformulated  in this new setting to a problem to describe the  morphism  from the winding quotient range  to the winding quotient of the defining  domain which is actually an elliptic curve for the moduli $q=\exp(2\pi\sqrt{-1}\tau)$. 

The new inversion problem asks to construct $q$-multiplicatively periodic functions on the winding quotient range, i.e.\ it asks to reformulate the classical elliptic function theory from the view point of the  winding quotient range. 
We solve the new  problem  by introducing  $q$-multiplicatively quasi-periodic functions $\mathcal{U}(w,q), \mathcal{V}(w,q)$ and a $q$-multiplicatively periodic function $\mathcal{W}(w,q)$ on the winding quotient range  whose pull-backs to the monodromy covering space are  the Fourier expansion parts of Weierstrass $\sigma$-function, $\zeta$-function and $\wp$-function, respectively,  up to a correction terms studied by Borcherds \cite{B}  generated by Eisenstein series $E_2$.
It seems interesting to note that this answer gives a connection between the classical elliptic function theory with the $q$-shifted factorial function $(w;q)_\infty$  in the study of combinatorics and special functions (see {\bf Appendix}).
Two areas of mathematics started from Euler meet again after 3\! centuries. 


\begin{rmk}
Actually, the pull-back of the function $\mathcal{W}(w,q)$ to the winding covering space has been studied in quite different contexts in mathematical physics  either as  the propagator $P(z)$ by Dijkgraaf \cite{D} (4.44) on elliptic curves and
the propagator $P_0^\infty(w_1,w_2;\tau,\overline{\tau})$ of BCOV theory by Si Li (\cite{L1} Lemma 3.1, \cite{L2} section 4 and \cite{LZ} 3.3), or as the genus-one Bergman kernel function $B(z_1,z_2)$  in the topological recursion by Eynard \cite{Ey} p.339 Examples (c.f.\  also Tyurin \cite{Ty} p.185 Example.\ $\omega_\alpha$ and Fay \cite{F} (46)).  
We do not know yet whether such connection of the winding phenomenon with mathmatical physics is common for non-geometric virtual period maps in general or it remains only in these test examples. We need to test further cases \cite{S3,G-S}.
\end{rmk}

\noindent
{\it Notation.}  (1)  We denote the set $\mathbb{C}\setminus\{0\}$ by $\mathbb{C}^\times$ because of its multiplicative structure. 

\noindent
(2)  When the complex plane $\mathbb{C}$, the punctured complex plane $\mathbb{C}^\times$, the complex upper half plane $\mathbb{H}$ or the unit disc $\mathbb{D}$ is equipped with a particular coordinate,  say $z$, then, in order to distinguish it from the other one equipped with other coordinates, we put to $\mathbb{C}$, $\mathbb{C}^\times$, $\mathbb{H}$ or $\mathbb{D}$ the coordinate name as a subscript, e.g.\ $\mathbb{C}_z$, $\mathbb{C}^\times_z$, $\mathbb{H}_z$ or $\mathbb{D}_z$.

\section{The virtual period map  $z^{\tau}$ on $\mathbb{C}_z^\times$}

We start to define  a virtual period map of rank 1 by the solution of period equation.\footnote{We use in \cite{S3} the term ``virtual'' for the period integrals whose associated hyperbolic root lattices do not have algebraic geometry origin, where the exponents become imaginary numbers. We employ in the present note the same term ``virtual'' only in the sense that {\it  the exponent $\tau$  of the period equation are non-real complex numbers}.}
Consider an ordinary differential equation in one variable $z$ and the unknown variable $w$:  
$$
  z\  \frac{dw}{dz} \ = \ \tau \ w 
$$
where $\tau$ is a complex number such that  $\rm{Im}(\tau)>0$. We  call the equation the {\it virtual period equation} of rank 1.\footnote{
If $\tau \in \frac{1}{2}\mathbb{Z}_{\ge0}$, the equation is the period  equation of type $\mathrm{A}_1$. }
 The solutions of the equation are, up to  a constant factor, equal to the function $w= z^{\tau}=\exp\!\big(\tau\log(z)\big)$. We say that the solution $z^\tau$  ``defines'', what we call,  the {\it virtual period map}:\!\!\!
$$
z^{\tau} \ : \  \mathbb{C}_z^\times \ \longrightarrow \ \mathbb{C}_{w}^\times,\ \  z \ \mapsto \ w=z^{\tau}
$$
from the punctured $z$-plane $\mathbb{C}_z^\times$ as the domain of the map  to  the punctured plane $\mathbb{C}_{w}^\times$ as the range of the map. Obviously it is not a single-valued map. So we define the map geometrically as follows.
We first choose $1\in \mathbb{C}_z^\times$ as the base point in the domain, and define the value of the map 
$z^{\tau}$ at $1\in  \mathbb{C}_z^\times$ to be 1.  Then, the values of the map are defined by analytic continuations of  $\exp\!\big(\tau\log(z)\big)$ along paths on the $z$-plane $\mathbb{C}_z^\times$ starting at 1.  The monodromy of this map by the counter clock-wise turn around $0$ of the $z$-variable is given by the multiplication of $\exp{\big(2\pi\sqrt{-1}\tau\big)}$ on $z^\tau$, i.e.\ its multiplication in the $w$-plane $\mathbb{C}_{w}^\times$. Since $|\exp{\big(2\pi\sqrt{-1}\tau\big)}|\!<\!1$, the monodromy group $\big(\exp{(2\pi\sqrt{-1}\tau)}\big)^\mathbb{Z}$ is an infinite cyclic group and gives a faithful representation of the fundamental group $\pi_1(\mathbb{C}_z^\times,1) \simeq \mathbb{Z}$. Then the monodromy covering space of the defining domain of the map $z^\tau$  is the universal covering space $\mathbb{C}_{\mathrm{Log}z}$ of $\mathbb{C}_z^\times$, where we denote by $\mathrm{Log}z$ the coordinate of the universal covering space such that $\exp(\mathrm{Log}z)=z$. 
On the space $\mathbb{C}_{\rm{Log}z}$, we have the covering transformation group action
 $$
  \Gamma_g\ :=\ 2\pi\sqrt{-1}\mathbb{Z}:  \qquad \mathrm{Log}z \mapsto \mathrm{Log}z +  2\pi\sqrt{-1}n \quad\ \ for \ \ n\in\mathbb{Z} \qquad \qquad  
  $$ 
 whose quotient map is the covering map  $\exp: \mathbb{C}_{\mathrm{Log}z}\to  \mathbb{C}_z^\times$.

\section{Lifting of the map $z^{\tau}$ to the monodromy covering space $\mathbb{C}_{\mathrm{Log}z}$}

Recall that $\mathrm{Log}z$ is the coordinate of the monodromy covering plane of $\mathbb{C}_z^\times$.
Then,  the virtual period map $z^{\tau}$ is lifted as the single-valued map $w=\exp\big(\tau\mathrm{Log}z\big)$ on the $\mathrm{Log}z$-plane $ \mathbb{C}_{\mathrm{Log}z}$.  So the range of virtual period map $z^{\tau}$ means the range of this lifted map. Thus, we obtain the diagram
$$\begin{array}{cccccc}
\vspace{0.2cm}
     \qquad  \mathrm{Log}z\  \in \!&  \mathbb{C}_{\mathrm{Log}z}\! &  \overset{ \exp(\tau \mathrm{Log}z) }{-\!\!\!-\!\!\!-\!\!\!-\!\!\!\longrightarrow}  & \mathbb{C}_{w}^\times\  \ni \  w=z^{\tau}=\exp(\tau{\mathrm Log}z) \\
     \vspace{0.2cm}
       \exp(\mathrm{Log}z) \!\!\!\!\!\!\!\!\!\!&  \downarrow  & \quad  \nearrow_{ z^{\tau}} \!\!\!\!\!\!\!\!\!\! \\
      \qquad   z\ \in\!\! \!\!\!\!\!\!\!&\!  \mathbb{C}_z\setminus\{0\}\!\!\!\!\!\!\!\!\!\!
\end{array}
$$
We warn that the anti-diagonal arrow for the virtual period map  $ z^{\tau}$ is a multi-valent map so that the diagram is not commutative in the usual sense. 

  The naive ``inversion problem'' in this setting is to give a map from $ \mathbb{C}_{w}^\times$ to $ \mathbb{C}_{z}^\times$ making the diagram commutative. However, {\it there does not exists a {\bf single-valued holomorphic map} 
which makes the diagram commutative}, since the analytic continuation of any local inversion map becomes  automatically the multi-valued function $z=w^{1/\tau}$ on $ \mathbb{C}_{w}^\times$. That is, {\it a naive solution of the classical inversion problem for the virtual period map $z^{\tau}$ by a single-valued map from the range to the domain does not exist.}

\section{Winding of the monodromy covering space $\mathbb{C}_{\mathrm{Log}z}$  to the range $\mathbb{C}_w^\times$
}
We saw that the virtual period map $z^{\tau}$ is presented by the single-valued function $w=\exp\big(\tau \mathrm{Log}z\big)$ defined on the monodromy covering space $\mathbb{C}_{\mathrm{Log}z}$ to the range  $\mathbb{C}_w$. Then, we immediately observe that its range $\mathbb{C}_w$ is the quotient of the monodromy covering space $\mathbb{C}_{\mathrm{Log}z}$ by the {\it translation group action}  
$$
\Gamma_w \ := \ \frac{2\pi\sqrt{-1}}{\tau}\mathbb{Z}: \qquad  \  \mathrm{Log}z \mapsto \mathrm{Log}z + \frac{2\pi\sqrt{-1}}{\tau}n \quad \ for \ \ n\in \mathbb{Z}  \qquad \qquad 
$$
%
%
where we note $\frac{2\pi\sqrt{-1}}{\tau} \in \mathbb{R}_{>0}$.
Let us call the group  $\Gamma_w$ the {\it winding group} and the map $\mathbb{C}_{\mathrm{Log}z}\to  \mathbb{C}_w^\times$ the {\it winding quotient map}.  
In the other words,  the range $\mathbb{C}_w^\times$ of $z^{\tau}$ is the {\it winding quotient space} of the monodromy covering space $\mathbb{C}_{\mathrm{Log}z}$\! by the winding group $\Gamma_w$-action.

 The key observation is that {\it  the group $\Gamma_w$ is caused  by the group  of the ``symmetry'' of the virtual period map $z^{\tau}=\exp(\tau\mathrm{Log}z)$.} Namely, the winding group action leaves the function $z^\tau$ invariant:
 $$
   \Big(e^{\frac{2\pi\sqrt{-1}}{\tau} n} z\Big)^{\tau}=z^{\tau} \quad \text{for } n\in \mathbb{Z} 
$$
or equivalently
$$
\exp\Big(\tau(\mathrm{Log}z+\frac{2\pi\sqrt{-1}}{\tau} n)\Big)=\exp\Big(\tau\mathrm{Log}z\Big) \quad \text{for } n\in \mathbb{Z}.
 $$
In the other words, $\Gamma_w$ does not change the branches of $z^\tau$, whereas $\Gamma_g$ does change them.


\section{Winding of the domain $ \mathbb{C}_z^\times$ of $z^{\tau}$ to the elliptic curve $E_\tau$}
In previous sections, we observed that the monodromy covering space $\mathbb{C}_{\mathrm{Log}z}$ obtained two different group actions, one: winding group $\Gamma_w$ and the other: the monodromy  covering group $\Gamma_g$. 
We observe immediately that the two actions commute to each other.  Consequently,  the winding group $ \Gamma_w$ acts also on the  defining domain $\mathbb{C}_z^\times $ of the map $z^{\tau}$ and the monodromy covering group $\Gamma_g$ acts on the range $\mathbb{C}_w^\times$.  Thus we obtain the following commutative Cartesian diagram: 
$$
(*)  \qquad\qquad\qquad
\begin{array}{cccccc}
   \vspace{0.05cm}      
  \qquad   &  \mathbb{C}_{\mathrm{Log}z}  &\!\!\!\!\!\!\!\!\!\!\!\!\!\!\!\!  \overset{ \bmod \Gamma_w }{-\!\!\!-\!\!\-\!\!\!-\!\!\!-\!\!\!-\!\!\!-\!\!\!-\!\!\!-\!\!\!\longrightarrow} \!\!\!\!\!\!\! & \!\!\!\!\!\!\!\!\!\!\!\!\!\!\!\!\!\!\!\!\!\!\!\!\!\!\!\!\!\!\! \mathbb{C}_w^\times \!\!\!\!\!\!\!\!\!\!\!\!\!\!\!\!\\    
          &\qquad \quad  \overset{\mid}{\underset{\downarrow}{\mid}} \bmod \ \Gamma_g &\   \nearrow_{  z^{\tau}}\!  &  \quad \overset{\mid}{\underset{\downarrow}{\mid}}  \bmod^\times \ \exp\big(\tau \Gamma_g\big)\!\!\!\!\!\!\!\!\!\!\!\!\!\!\!\!\!\!\\
          
         &\!\!  \mathbb{C}_z^\times  &\!\!\!\!\!\!\!\!\!\!\!\!\!\! \overset{ \bmod^\times  \exp(\Gamma_w) }{-\!\!\!-\!\!\!-\!\!\-\!\!\!-\!\!\!-\!\!\!-\!\!\!-\!\!\!-\!\!\!-\!\!\!\longrightarrow} \!\!\!\!\! & \!\!\!\!\!\!\!\!\!\!\!\!\!\!\!\!\!\!\!\!\!\!\!\!\!\!\!\!\!\!\! E_\tau \!\!\!\!\!\!\!\!\!\!\!\!\!\!\!\!\!\!\\
\end{array}
\qquad\qquad\qquad\qquad\qquad\qquad\qquad\qquad\qquad
$$
where   $E_\tau$ is the quotient of the domain $  \mathbb{C}_z^\times$ by the multiplicative winding group $\Gamma_w$-action. Actually, $E_\tau$ is an elliptic curve given by the  quotient of the plane $\mathbb{C}_{\mathrm{Log}z}$  by the action of the lattice given by 
$$
\Gamma_g+\Gamma_w\ = \  \frac{2\pi\sqrt{-1}}{\tau}(\mathbb{Z}+\mathbb{Z}\tau).
$$ 


The virtual period map $z^\tau$ induces the new virtual period  map  $: E_\tau \to \mathbb{C}_w^\times$ (recall that $z^\tau$ is invariant under the action of $\Gamma_w$) so that $\mathbb{C}_w^\times$  becomes the monodromy covering space of $E_\tau$ for the new virtual period map,  since $\pi_1(E_\tau,0)\simeq \Gamma_g+\Gamma_w$ and the kernel of the new monodromy presentation is equal to $\Gamma_w$.


\smallskip
In this new situation for the virtual period map, what does the ``inversion problem'' may mean? One answer, which we shall employ in the present note, is that it is a question asking to describe the inverse morphism $\mathbb{C}_w^\times \to E_\tau$, i.e.\ the vertical down arrow in the RHS of the diagram $(*)$.  Namely, for a prescribed $w\in \mathbb{C}_w^\times$, we want to find the {\it unique}  point, say $e$, in  $E_\tau$ whose inverse image  in $\mathbb{C}_w^\times$ (=the set of virtual periods for  $e$) contains $w$, and to describe $e$ in terms of $w$. 

 In order to realize this program, we need to fix a coordinate presentation\! of the point $ e \in E_\tau$.  It is well-known that the elliptic curve $E_\tau$ is realized as the plane curve:\!\!\!\!
$$
\qquad \qquad E_\tau \ \subset \  \mathbb{P}^2, \quad \  u \bmod \mathbb{Z}+\mathbb{Z}\tau \ \ \mapsto \ \ (\wp(u,\tau):\wp'(u,\tau):1)
$$ 
by the use  of Weierstrass $\wp$-function $\wp(u,\tau)$ and its derivative $\wp'(u,\tau)$ given in  the coordinate 
$$
u=\frac{\tau}{2\pi\sqrt{-1}}\mathrm{Log}z
$$ 
where the scaling factor $\frac{\tau}{2\pi\sqrt{-1}}$ comes from the scaling of the group $\Gamma_g+\Gamma_w= \frac{2\pi\sqrt{-1}}{\tau}(\mathbb{Z}+\mathbb{Z}\tau)$ relative to the lattice $\mathbb{Z}+\mathbb{Z}\tau$ for the $\wp$-function.\footnote
{ The realization of $E_\tau$ as a plane curve is not unique (e.g.\ \cite{S2}). The present note has chosen  $\mathrm{A}_2$-type realization. It would be interesting to study the other cases $\mathrm{B}_2$-type and $\mathrm{G}_2$-type.
}
Thus, {\it the new inversion problem asks to describe the pair coordinates  $(\wp(u,\tau),\wp'(u,\tau))$ of the point $e \in E_\tau$ in terms of $w=\exp(\tau\mathrm{Log}z)=\exp(2\pi\sqrt{-1} u)$ and $\tau$.}

\section{$q$-multiplicatively quasi-periodic function $\mathcal{V}(w,q)$ and  periodic function $\mathcal{W}(w,q)$}

In order  to answer to the new inversion problem, we review the classical elliptic function theory from the view point of $q$-periodic, i.e.\ invariant by the multiplicative action of $q^{\mathbb{Z}}\!:=\!\exp(\tau\Gamma_g)\!=\!\exp(2\pi\sqrt{\!-1} \tau\mathbb{Z})$,  functions on the winding quotient space $\mathbb{C}_w\times\mathbb{C}_q^\times$. We start with  
$q$-shifted factorial function $(w,q)_\infty$, which we regard as the holomorphic function on the winding quotient space. Its logarithmic derivative induces a ``positive'' $q$-multiplicatively quasi-periodic function $\mathcal{V}_+(w,q)$. By making up the ``negative'' part $\mathcal{V}_-(w,q)$ together, we obtain the
$q$-quasi-periodic  function $\mathcal{V}(w,q)$ of the period 1, which, as a meromorphic function, has simple poles along  $1\!-\!wq^m\!=\!0$ for $m\! \in\! \mathbb{Z}$. The logarithmic derivative  $\mathcal{W}(w,q)$ of $\mathcal{V}(w,q)$ is a $q$-multiplicatively periodic meromorphic function on $\mathbb{C}_w^\times\times \mathbb{D}_q^\times$, which gives the answer to the new inversion problem.\
 See {\bf Appendix} for more relationships of them  with elliptic function theory. 

We start with an infinite product (see Euler \cite{E} p.256)  in  two variables $w$ and $q$
\vspace{-0.15cm}
$$
(1) \qquad \qquad\qquad\qquad \qquad\qquad\quad\quad  
(w;q)_\infty \ =\ \prod_{m=0}^\infty (1-wq^m)
\qquad\qquad\qquad \qquad\qquad\qquad \qquad  \qquad\ \ \
\vspace{-0.1cm}
$$
called  $q$-shifted factorial function (\cite{Sl} p.88, \cite{GR} p.351).  Let us show that {\it the infinite product {\rm (1)}  converges absolutely and compact uniformly  to a holomorphic function on  the winding quotient space $\mathbb{C}_w\times \mathbb{D}_q$, whose zero-loci is the union of the divisor $\cup_{n=0}^\infty\{1-wq^n=0\}$. We shall denote the convergent function by the same  $(w,q)_\infty$.}
\begin{proof} 
%
For any compact set $K\subset \mathbb{C}_w\times \mathbb{D}_q$, set $c:=\sup\{|w| \mid  (w,q) \in K\}$ and $d\!:=\!\sup\{|q| \! \mid \! (w,q)\!\in\! K\}$ so that $\sup\{ |wq^n| \! \mid \!  (w,q)\!\in\! K\}\! \le \! cd^n$.  Since $0\!\le \! d\! <\!1$, this bound is sumable, i.e.  $\sum_{n=0}^\infty cd^n \!=\! c/(1\!-\!d)<\infty$. This implies the absolute uniform convergence of (1) to a holomorphic function on $K$ 
(see  \cite{WW} p.30 for infinite products).
\end{proof}

As a consequence, we introduce three infinite product functions using $(w,q)_\infty$.\footnote{ Since $\mathcal{U}_-(w,q) $ and $\mathcal{U}(w,q) $ (and the following  $\mathcal{V}_-(w,q) $, $\mathcal{V}(w,q) $ and $\mathcal{W}(w,q) $)  have essential singularities at $w=0$, therefore, from now on, we restrict their defining domain to $\mathbb{C}_w^\times\times\mathbb{D}_q^\times$.
}
\vspace{-0.08cm}
$$
\begin{array} {ccccl}
\!\qquad\qquad   \mathcal{U}_+(w,q) & :=& (w;q)_\infty &\! =\!& \prod_{m=0}^\infty (1-wq^m).  \qquad\qquad\qquad \qquad\qquad\qquad \qquad\qquad\qquad\quad \\
\!(2) \qquad\quad \!  \mathcal{U}_-(w,q) & :=& (w^{-1};q)_\infty/(1-w^{-1}) &\! =\!& \prod_{m=1}^\infty (1-w^{-1}q^m).  \qquad\qquad\qquad \qquad\qquad\qquad \qquad\qquad\qquad\qquad \\
  \!  \qquad\qquad   \mathcal{U}(w,q) & :=\! &\!\! \mathcal{U}_+(w,q)\ \mathcal{U}_-(w,q) &\! =\! & \prod_{m=0}^\infty (1\!-\!wq^m) \prod_{m=1}^\infty (1\!-\!w^{-1}q^m).  \qquad\qquad\qquad \qquad\qquad\qquad \qquad\qquad\qquad
\end{array}
$$
Since they converge compact uniformly to  holomorphic functions on {\small  $\mathbb{C}_w^\times\times \mathbb{D}_q^\times$}, we can derivate them by the variable $w$ term-wisely. Thus we obtain three  infinite sums 
\vspace{-0.1cm}
$$
\begin{array}{clccl}
 \qquad\quad   \mathcal{V}_+(w,q) &:= & - w\frac{\partial}{\partial w} \log(\mathcal{U}_+(w,q)) & =&    \sum_{m\in\mathbb{Z}_{\ge0}}\frac{wq^m}{1-wq^m} . \qquad\qquad \qquad\qquad\qquad \qquad\qquad\qquad\quad \\ 
\!(3)  \qquad   \mathcal{V}_-(w,q) &:= & -w\frac{\partial}{\partial w} \log(\mathcal{U}_-(w,q)) & =&    \sum_{m\in\mathbb{Z}_{<0}}\frac{1}{1-wq^m} . \qquad\qquad \qquad\qquad\qquad \qquad\qquad\qquad\quad \\ 
 \qquad\quad   \mathcal{V}(w,q) &:= & -w\frac{\partial}{\partial w} \log(\mathcal{U}(w,q)) & =&    \sum_{m\in\mathbb{Z}_{\ge0}}\frac{wq^m}{1-wq^m} +    \sum_{m\in\mathbb{Z}_{<0}}\frac{1}{1-wq^m}. \qquad\qquad \qquad\qquad\qquad \qquad\qquad\qquad\quad \\ 
\end{array}
\vspace{-0.08cm}
$$
Note that, by definition directly,  we have following formal relations
\vspace{-0.08cm}
$$
(4) \qquad \qquad \quad \mathcal{V}_-(w,q)=- \mathcal{V}_+(w^{-1},q)- \frac{1}{1-w}  \quad \text{and}\quad 
\mathcal{V}(w,q) = \mathcal{V}_+(w,q) +\mathcal{V}_-(w,q). \qquad  \qquad
\vspace{-0.08cm}
$$
Also by the definition, we obtain following involution property of  $\mathcal{U}(w,q)$ and $\mathcal{V}(w,q)$:
$$
(5)  \quad \qquad \qquad   \qquad \mathcal{U}(w,q)/\mathcal{U}(w^{-1},q)=-w  \quad \quad    \mathcal{V}(w,q)\ + \ \mathcal{V}(w^{-1},q)\ + \ 1 \ =\ 0.    \quad  \quad \quad  \quad \quad   \qquad   
$$

The absolute compact uniform convergence of (1) implies that these functions are holomorphic except along the singular loci: the  divisor $\cup_{n=0}^\infty\{1-wq^n=0\}$.
In order to confirm that they are meromorphic functions on  $\mathbb{C}_w^\times\times \mathbb{D}_q^\times$ along the singular loci, we show the compact uniform convergence of the sums (2) in a generalized sense.\!\!\!


\bigskip
\noindent
{\bf Assertion A.} 
 (i)  {\it The sums {\rm (3)} converges compact uniformly, in a generalized sense given in the proof, to the  meromorphic functions on  $\mathbb{C}_w^\times\times \mathbb{D}_q^\times$. }

\smallskip
\noindent


\noindent
(ii)  {\it The functions {\rm (2)} and {\rm (3)} are $q$-quasi periodic in the following sense:}
$$
\begin{array}{cclccc}
 \quad\qquad \mathcal{U}_+(w,q)&=& \mathcal{U}_+(qw,q)(1-w)&  \mathcal{V}_+(w,q)  &\!\!\!\! \!\!=\!\!\!\!& \!\!\!\!\!  \mathcal{V}_+(qw,q)\ \ + \  \ \frac{w}{1-w}
\qquad\qquad\qquad\qquad\qquad\qquad\qquad \quad \\
\!\!\!(6)  \quad\quad \mathcal{U}_-(w,q)&=& \mathcal{U}_-(qw,q)(1-w^{-1})^{-1}&\quad  \mathcal{V}_-(w,q) \quad&\!\! \!\!\!\!=\!\!\!\!& \quad \ \mathcal{V}_-(qw,q)\ \ - \ \ \frac{1}{1-w}
\qquad\qquad\qquad\qquad\qquad \qquad\qquad\qquad \quad \\
 \quad\qquad \mathcal{U}(w,q)&=& \ \mathcal{U}(qw,q)(-w)&  \mathcal{V}(w,q)   &\!\!\!\!\!\!= \!\!\!\!& \!\! \mathcal{V}(qw,q)\ \ \ - \ \ \ \   1
\qquad\qquad\qquad\qquad\qquad\qquad\qquad \quad \ \ 
\end{array}
 $$
 \noindent
\begin{proof}
\!\!(i)  We first show only the case $ \mathcal{V}_+(w,q)$ of (3). That is, in the following proof, (3) means only the case $ \mathcal{V}_+(w,q)$.  
First, let us fix following  notations.\!\!
 
 \noindent
 (a) For $n\in \mathbb{Z}_{\ge0}$, set $\mathcal{U}_n:=\{(w,q)\in \mathbb{C}_w\times \mathbb{D}_q^\times\mid  |w|<|q|^{-n-1}\}$.  
 
 \noindent
 \quad\ \  Then  $\mathcal{U}_n$ ($n\in\mathbb{Z}_{\ge0}$) is an increasing sequence such that $\mathbb{C}_w^\times\times \mathbb{D}_q^\times=\cup_{n\in \mathbb{Z}_{\ge0}} \mathcal{U}_n$. 
 
 \noindent
 (b) For any compact subset $K$ of the domain $\mathcal{U}_n$ and for any complex valued 

\noindent
\quad \ \ 
 continuous function $f$ on $\mathcal{U}_n$, we set $|\!|f|\!|_K:=\mathrm{sup}\{|f(w,q)|\mid (w,q)\in K\}$.

 \medskip
 \noindent
Next, we show that the sum (3)$_n$, obtained from (3) by removing the first $n+1$-terms for $m\!=\! 0,1,\cdots, n$,  converges uniformly on any compact subset $K$ of the  domain $\mathcal{U}_n$. 

\smallskip
{\it Proof.}  
Since $|\!|q|\!|_K<1$, we have  $|\!|wq^m|\!|_K \le |\!|wq^{n+1}|\!|_K<1$ for $m>n$.  Then  $|\!|\frac{1}{1-wq^m}|\!|_K\le \frac{1}{1- |\!|wq^m|\!|_K}\le  \frac{1}{1- |\!|wq^{n+1}|\!|_K}$. 
Thus, we have the following majoration 
$$
\sum_{m\in\mathbb{Z}_{>n}}   |\!| \frac{wq^m}{1-wq^m}|\!|_K  \ \le  \ \frac{ |\!|wq^n |\!|_K}{1- |\!|wq^{n+1}|\!|_K}\sum_{m>n} |\!|q|\!|_K^{m-n} \ = \  
  \frac{ |\!|wq^n |\!|_K}{1- |\!|wq^{n+1}|\!|_K} \frac{|\!|q|\!|_K}{1-|\!|q|\!|_K}.
\vspace{-0.1cm}
$$
which implies the compact uniform convergence of the sum (3)$_n$. \qquad \qquad \qquad $\Box$

 Thus, (3)$_n$ gives a holomorphic function on $\mathcal{U}_n$. Then the total sum (3) defines a meromorphic function on $\mathcal{U}_n$, having the simple poles along $1\!-\!wq^m\!=\!0$ ($m=0,1,\cdots,n$). 
Recalling that the sequence $\mathcal{U}_n$ ($n\in \mathbb{Z}_{\ge0}$)  exhausts  $\mathbb{C}_w\times \mathbb{D}_q^\times$, we conclude that the sum (3) defines a meromorphic function on $\mathbb{C}_w\times \mathbb{D}_q$.  We say that  (3) converges compact uniformly to the meromorphic function $\mathcal{V}_+(w,q)$ in a generalized sense.

The convergence of two other  cases  $\mathcal{V}_-(w,q)$ and $\mathcal{V}(w,q)$ of (3)  follows from (4). 

\medskip
\noindent
{\it Note.}  By the construction, the principal parts of poles of  $ \mathcal{V}_+(w,q)$  are given by $\frac{1}{1-wq^m}$ ($m\in \mathbb{Z}_{\ge0}$) whose residues are $-q^{-m}$. 
 
 
(ii)    The quasi-periodicity (6) of the function $\mathcal{V}_+(w,q)$  follows, since the $q$-shift $\mathcal{V}_+(wq,q)$ is actually expressed by the sum (3)$_0$ and the 0-th term of (3) is  $\frac{w}{1-w}$. The quasi-periodicity  of two other  cases  $\mathcal{V}_-(w,q)$ and $\mathcal{V}(w,q)$ follows from (4). 
\end{proof}

 \begin{rmk}
  Naive candidates of the functions $\mathcal{U}(w,q)$ and  $\mathcal{V}(w,q)$ may be the  infinite product $\prod_{m\in \mathbb{Z}}(1-wq^m)$ and the infinite sum  $ \sum_{m\in\mathbb{Z}}\frac{wq^m}{1-wq^m}$, respectively, which are formlly $q$-periodic. However, their  negative parts of the index $m\in \mathbb{Z}_{<0}$ do not converge. Therefore, we defined $\mathcal{U}(w,q)$ and $\mathcal{V}(w,q)$ by replacing  the negative parts by the convergent $\mathcal{U}_-(w,q)=  (w^{-1};q)_\infty/(1-w^{-1}) $ and  $\mathcal{V}_-(w,q)= \sum_{m\in\mathbb{Z}_{<0}}\frac{1}{1-wq^m} $, respectively, as in (2) and (3). These modifications break the (formal) $q$-periodicity of the series, but cause the quasi $q$-periodicity of the functions  $\mathcal{U}(w,q)$ and $\mathcal{V}(w,q)$ as in (6). We shall see these two ``problems'' (i.e.\ the non-convergency and the non-periodicity) are resolved in the  function $\mathcal{W}(w,q)$ introduced in the next paragraph.
 %
 \end{rmk}
 
 We now apply again  the action of  $w\frac{\partial}{\partial w}$ on $ \mathcal{V}(w,q) $. 
Since the sum (3) is compact uniform convergent (in the generalized sense given in the proof of {\bf Assertion A}), we can apply the derivation to the series term-wisely so that we obtain  a new series:
$$
\vspace{-0.2cm}
(7)  \qquad \qquad\qquad \qquad\qquad \qquad  \mathcal{W}(w,q)\  := \  w\frac{\partial}{\partial w} \mathcal{V}(w,q) \ =  \ \sum_{m\in\mathbb{Z}}\frac{wq^m}{(1-wq^m)^2}  . \qquad\qquad \qquad  \qquad\qquad
$$
which is compact uniformly convergent for all $m\in\mathbb{Z}$  in the generalized sense.
Noting the action $w\frac{\partial}{\partial w}$ is invariant by the $q$-multiplication $w\!  \mapsto \! qw$, the $q$-quasi-periodi-city (6) of $ \mathcal{V}(w,q)$ implies the {\it $q$-multiplicatively periodicity of the function $\mathcal{W}(w,q)$}. 
$$
(8)  \qquad\qquad\qquad\qquad   
\mathcal{W}(qw,q)\   = \  \mathcal{W}(w,q) \qquad \text{and} \qquad  \mathcal{W}(w,q)\   =   \ \mathcal{W}(w^{-1},q)  .\quad \qquad\qquad \qquad\qquad\qquad\qquad\qquad 
$$

\noindent
{\bf Assertion B.} {\it The function $\mathcal{W}$ on  $\mathbb{C}_{w}^\times\times \mathbb{D}_q^\times$ and the $\wp$-function on $\mathbb{C}_{u} \times \mathbb{H}_\tau$ are related by the following pull-back relations.
\vspace{-0.2cm}
$$
(9) \qquad\qquad\qquad \qquad \qquad  \ \ 
\wp(u,\tau) \ = \ (2\pi\sqrt{-1})^2\ \mathcal{W}(w,q) \ -\ \frac{\pi^2}{3}E_2(\tau) 
  \quad \quad  \quad\quad\qquad\qquad\qquad\quad 
\vspace{-0.1cm}
$$
where we set 
$$
(10)\qquad\qquad\qquad\qquad\qquad\quad\  w=\exp(2\pi\sqrt{-1}u), \qquad  q=\exp(2\pi\sqrt{-1}\tau), \qquad\qquad\qquad\quad \qquad\qquad
$$
and
 $E_2(\tau)=\frac{1}{\pi^2}\sum_m \big(\sum_n' \frac{1}{(m\tau+n)^2}\big)$ 
 is the Eisenstein series of weight 2 (here $\sum_n'$ means $n\!=\!0$ is omitted if $m\!=\!0$,  see} Zagier \cite{Z} p.\ 19).

\begin{proof} 
For each fixed $\tau\in \mathbb{H}$, the substitution of  $w=\exp(2\pi\sqrt{-1}u)$ and $ q=\exp(2\pi\sqrt{-1}\tau)$ to $\mathcal{W}(w, q)$ is a doubly periodic meromorphic function in $u\in \mathbb{C}_u$ having poles at the lattice $\mathbb{Z} + \mathbb{Z}\tau$, and the principal part of its Laurent expansion at $u=0$ is equal to $\frac{1}{(2\pi\sqrt{-1})^2}\frac{1}{u^2}$. Then, the difference of the $\wp(u,\tau)$ and the substitution $(2\pi\sqrt{-1})^2\mathcal{W}(\exp(2\pi\sqrt{-1}u), \exp(2\pi\sqrt{-1}\tau))$ is a doubly periodic function on $u\in \mathbb{C}_u$ without poles so that it is a constant function on $u$. That is, the difference is a holomorphic function depending only on $\tau$.   

Actually, the differnce was already known by R.\ Borcherds as follows.
In \cite{B} Chapter 7,  Borcherds give the following Fourier expansion of the $\wp$-function  (see also  Aoki \cite{A} (3.10) and Lang \cite{L} Chap.4 Sec.2 Prop.\ 2-3).
\vspace{-0.1cm}
$$
\wp(u,\tau)=(2\pi\sqrt{-1})^2\sum_{m\in\mathbb{Z}}\sum_{n\in \mathbb{Z}_{>0}}n \exp\big(\pm2\pi\sqrt{-1}n(u+m\tau)\big)  -\frac{\pi^2}{3}E_2(\tau)
$$
whenever $\Im(u+m\tau)$ is non-zero for all integers $m$, where 
 the sign $\pm$ means $+$ if  
$|\exp(2\pi\sqrt{-1}n(u+m\tau))|<1$ and $-$ if  
$|\exp(2\pi\sqrt{-1}n(u+m\tau))|>1$. 

That is, the pull-back relation holds on the open set $\mathbb{C}_u\times \mathbb{H}_\tau \setminus \cup_{m\in \mathbb{Z}} \{\Im(u+m\tau)=0\}$. 
Clearly, this set projects surjectively onto  $\mathbb{H}_\tau$, implying the equality holds.
\end{proof}   
 
 \begin{rmk} \ \  The $q$-multiplicative periodicity of $\mathcal{W}(w,q)$ comes from the monodromy group $\Gamma_g$ action on the covering space $\mathbb{C}_{\rm{Log}z}$ (whose invariant is the exponential function)  so that one should not confuse it with the winding group $\Gamma_w$ invariance. Symbolically speaking, {\it  $\mathcal{W}(w,q)$ is the winding quotient of the exponential map $\exp: \mathbb{C}_{Loz} \to \mathbb{C}_z^\times$, i.e.\ $\mathcal{W} = \exp \bmod  \ \Gamma_w$.} In the other words, the function $\mathcal{W}$ is the shadow of the exponential function by the winding group $\Gamma_w$.  
 \end{rmk}
 
  \begin{rmk}  \   Yosuke Ohyama has informed author the followings
  
 \noindent
  1.\  The function $\mathcal{U}(w,q)$  gives the double product part of Jacobi triple product:
  $$
  \sum_{m=-\infty}^{m=\infty} q^{\frac{m(m-1)}{2}} (-w)^m= \prod_{m=1}^\infty (1-q^m) 
  \prod_{m=0}^\infty (1-wq^m)  \prod_{m=1}^\infty(1-w^{-1}q^m)
  $$
    
 \noindent
  2.\ The\! function\! $\mathcal{V}_+(w,q)$\! is\! expressed\! by\! a\! Lambert\! series\! \cite{Kn}\! p.452 (f) on\! the\! disc\! $\mathbb{D}_w$.\!\!\!
  \vspace{-0.2cm}
  $$
 \quad \quad  \sum_{n=0}^{\infty} \frac{w^n q^n }{1-q^n}  \ = \sum_{n=0}^{\infty} \frac{w q^n}{1-w q^n} \Big|_{\mathbb{D}_w\times \mathbb{D}_q} \ =\ \mathcal{V}(w,q) \Big|_{\mathbb{D}_w\times \mathbb{D}_q} 
\vspace{-0.2cm}
 $$

\noindent
 3.\ The\! function\!  $\mathcal{W}(w,q)$\! is\! expressed\! by\!  bilateral\! basic\! hypergeometric\! series\! (\cite{GR}\! (chap.5.).\!\!\!
     $$
 {}_2\psi_2
 \begin{bmatrix}
 w ,w & \\
 \vspace{-0.2cm}
 &\!\!; q,q \\
 wq,wq\!\!&
 \end{bmatrix}
 =
\! \sum _{n=-\infty}^{\infty}\!
 \frac{(w,w;q)_n}{(wq,wq;q)_n}
 q^n
 =
 \! \sum _{n=-\infty}^{\infty}\!
 \frac{(1-w)^2}{(1-wq^n)^2}
 q^n
 =\frac{(1-w)^2}{w} \mathcal{W}(w,q).
 $$

\end{rmk}

\smallskip 
\section{Solution to the new inversion problem}
Using the $q$-multiplicatively periodic function $\mathcal{W}(w,q)$ constructed in the previous section, the descriptions of the morphism: $\mathbb{C}_w^\times \to E_\tau\subset \mathbb{P}^2$ is achieved as follows.  This  gives an answer to the new inversion problem posed at the end of  section 5.
\begin{thm}
For each fixed $q=\exp(2\pi\sqrt{-1}\tau)$, the pair of functions:
$$(\ (2\pi\sqrt{-1})^2\mathcal{W}(w,q) \!-\!\frac{\pi^2}{3}E_2(\tau)\ :\  (2\pi\sqrt{-1})^3w\frac{\partial}{\partial w}\mathcal{W}(w,q) \ :\ 1\ )
$$
 gives a morphism: $\mathbb{C}_w^\times \to E_\tau \subset \mathbb{P}^2$ making  the diagram $(*)$ in section 5 commutative.
%
\end{thm}
\begin{proof}  This follows from {\bf Assertion B.} in section 6.  \end{proof}

\section{ {\bf Appendix }}
We describe the relationship between Euler's symbol $(w;q)_\infty$  and  the Weierstrass elliptic functions.  We use the results and  the equations (1)\! -\! (10) in section 6. 

\medskip
Recall that, starting from Euler's symbol $(w;q)_\infty= \prod_{m=0}^\infty (1-wq^m)$  (1),  we have introduced the (either holomorphic or meromorphic) functions on $(w,q)\in \mathbb{C}_w^\times\times\mathbb{C}_q^\times$: 
$$
\begin{array}{cccccl}
\!\!\!\!\!\!(2) \qquad & \mathcal{U}(w,q) &\!=\!& \frac{(w;q)_\infty \cdot(w^{-1},q)_\infty}{(1-w)} &\!=\!& \prod_{m=0}^\infty (1\!-\!wq^m) \prod_{m=1}^\infty (1\!-\!w^{-1}q^m) \quad \\
\!\!\!\!\!\!(3)\qquad  &   \mathcal{V}(w,q)& \!=\!&  -w\frac{\partial}{\partial w} \log(\mathcal{U}(w,q))  &\!= \!&   \sum_{m\in\mathbb{Z}_{\ge0}}\frac{wq^m}{1-wq^m} +    \sum_{m\in\mathbb{Z}_{<0}}\frac{1}{1-wq^m} \\
\!\!\!\!\!\!  (7)\qquad  &  \mathcal{W}(w,q) &\!=\!&  w\frac{\partial}{\partial w} \mathcal{V}(w,q) &\!= \!&\sum_{m\in\mathbb{Z}}\frac{wq^m}{(1-wq^m)^2} 
   \end{array}
$$
  and have shown the relationship 
$$
(9) \qquad\qquad\qquad \qquad \qquad  \ \ 
\wp(u,\tau) \ = \ (2\pi\sqrt{-1})^2\ \mathcal{W}(w,q) \ -\ \frac{\pi^2}{3}E_2(\tau) 
  \quad \quad \quad  \quad \quad\qquad\qquad\qquad\qquad\qquad 
$$
where the coordinates $(w,q)$ for the function $\mathcal{U}(w,q),\mathcal{V}(w,q)$ and $\mathcal{W}(w,q)$ and the coordinates $(u,\tau)$ for the $\wp$-functions are related by the winding quotient morphism
$$
(10) \qquad\qquad\qquad\qquad\qquad\qquad\  w=\exp(2\pi\sqrt{-1}u), \qquad  q=\exp(2\pi\sqrt{-1}\tau),  \qquad\qquad\qquad\qquad\qquad 
$$

It is straight forward now to  show the following.

\bigskip
\noindent
{\bf Assertion C.} {\it Weierstrass $\zeta$-function and $\sigma$-function are described by the functions $\mathcal{V}(w,q)$ and $\mathcal{U}(w,q)$ as follows, where we use the coordinate change} (10).\!\!\!\! 
$$
(11) \qquad \qquad \qquad \qquad
\zeta(u,\tau) \ = \ -2\pi\sqrt{-1}\ \mathcal{V}(w, q) \ + \  \frac{\pi^2}{3}E_2(\tau) u  \ - \ \pi\sqrt{-1}  \qquad  \qquad
 \qquad \quad \qquad  \quad \ \ \ 
 $$
 \vspace{-0.6cm}
$$
(12) \qquad \qquad \qquad  \qquad   \sigma(u,\tau) \ = -\frac{1}{2\pi\sqrt{-1}} \frac{ \mathcal{U}(w,q)}{\eta(\tau)^2} \exp\Big(\frac{\pi^2}{6} E_2(\tau)u^2 - \pi\sqrt{-1}u+\frac{\pi\sqrt{-1}}{6}\tau\Big)  .
\quad \quad \qquad    \qquad 
$$
 {\it where $\eta(\tau)$ is Dedekind eta-function: $\exp(\frac{\pi\sqrt{-1}}{12}\tau)\prod_{m=1}^\infty(1-q^m)$.}

\begin{proof} 
Recall that Weierstrass $\zeta$-function and $\sigma$-function are defined by the equations   $\wp(u,\tau)= -\frac{\partial}{\partial u} \zeta(u, \tau)$ and $\frac{\partial}{\partial u} \log \sigma(u,\tau)=\zeta(u,\tau)$  together with some initial conditions. Note also that  we have $\frac{\partial}{\partial u}=2\pi\sqrt{-1} w\frac{\partial}{\partial w}$ due to  (10). 

\medskip
Case for $\zeta(u,\tau)$.

\noindent
In view of (7) and (9), we have
$
\frac{\partial}{\partial u}  \zeta(u,\tau) +(2\pi\sqrt{-1})^2w\frac{\partial}{\partial w} \mathcal{V}(w,q) - \frac{\pi^2}{3}E_2(\tau)=0
$.
We integrate it by $du\!=\!\frac{1}{2\pi\sqrt{-1}}\frac{dw}{w}$, where  we set the integral constant to be $F(\tau)$. That is,\!\!\! 
\vspace{-0.2cm}
$$
(*) \qquad \qquad\qquad\qquad  \zeta(u,\tau) \ + \ 2\pi\sqrt{-1} \mathcal{V}(w, q) \ - \  \frac{\pi^2}{3}E_2(\tau) u  \ =\ F(\tau).  \qquad \qquad \qquad
\vspace{-0.1cm}
$$ 
Recall that the zeta function is an odd function w.r.t.\  the variable $u$.  On the other hand, the involution  $u\leftarrow \!\!\! \rightarrow -u$ induces the involution $w \leftarrow\!\!\! \rightarrow w^{-1}$.  We apply the involutions to $(*)$. Recalling the involution relation (5) on $\mathcal{V}(w,q)$, we obtain  
\vspace{-0.1cm}
$$
\quad -\zeta(u,\tau)\  -\  2\pi\sqrt{-1} (\mathcal{V}(w, q) \ + \ 1)\ + \ \frac{\pi^2}{3}E_2(\tau) u 
=  F(\tau)
\vspace{-0.1cm}
$$
Suming this with $(*)$, we get  $2F(\tau)= -2\pi\sqrt{-1}$, and the formula (11).

\medskip
Case for $\sigma(w,\tau)$.

\noindent
In view of (3) and (11), we have
$      
\frac{\partial}{\partial u} \log \sigma(u,\tau) =2\pi\sqrt{-1}\ w\frac{\partial}{\partial w} \log \mathcal{U}(w, q) +   \frac{\pi^2}{3}E_2(\tau) u   -  \pi\sqrt{-1}
 $.
 We integrate it again by $du=\frac{1}{2\pi\sqrt{-1}}\frac{dw}{w}$, where we set  the multiplicative integral constant to be $G(\tau)$. That is, 
 \vspace{-0.1cm}
$$
 \qquad \qquad \qquad \qquad \ \   \sigma(u,\tau) = \mathcal{U}(w,q) \exp\big(\frac{\pi}{6} E_2(\tau)u^2 - \pi\sqrt{-1}u\big)G(\tau)
\qquad \qquad \qquad\ \    
\vspace{-0.1cm}
$$
Recall that $\underset{u\to0}{\lim}\frac{\sigma(u,q)}{u}=1$, where ``$u\to 0$'' is the same as ``$w\to 1$''. Thus, we obtain 
\vspace{-0.3cm}
$$
1= \underset{w\to 1}{\lim}\frac{2\pi\sqrt{-1}\mathcal{U}(w,q)}{\log(w)}  \cdot G(\tau)
\vspace{-0.1cm}
$$
where we take the principal value of $\log(w)$.  Recall the Tayler expansion $\log(w)= (w-1) -   \frac{1}{2}(w-1)^2 + \cdots$ at $w=1$. Using (2), we obtain
$$
2\pi\sqrt{-1} \frac{\mathcal{U}(w,1)}{w-1}\Big|_{w=1}=-2\pi\sqrt{-1}\prod_{m>0}(1-q^m)^2= - 2\pi\sqrt{-1}\exp(-\pi\sqrt{-1}\tau/6)\eta(\tau)^2,
\vspace{-0.3cm}
$$
so that $G(\tau)=-\frac{\exp(\frac{\pi\sqrt{-1}}{6}\tau )}{2\pi\sqrt{-1}\eta(\tau)^2}$. 
So, we get the formula (12). 

 \end{proof}
 
 \noindent
 {\it Note.} That $\sigma(w,\tau)$ is an odd function in $u$ follows immediately from (12) and (5).  However it\! does\! not\! determine the constant term $\frac{\pi\sqrt{-1}}{6} \tau$ in the exponential part of (12).\!\!\!
 
 \smallskip

\medskip
\noindent
{\it Acknowledgements.}
The author express his gratitude to Hiroki Aoki who informed him the Borcherds work (see the  proof of {\bf Assertion B}) and helped him in solving the new inversion problem,  to Akishi Ikeda who explained him the roles of the pull-back of the function $\mathcal{W}$ to the winding covering space in mathematical physics and informed him its literatures (see Remark 1.1), and to Yosuke Ohyama for guiding to bilateral hypergeometric functions and related interesting discussions on references and histories (see Remark 6.3). He also express his gratitude to Masahiko Yoshinaga for discussions during the preparation of the present note.



\end{document}